\documentclass{elsart}
\usepackage{amssymb}
\usepackage{graphicx}

\begin{document}

\begin{frontmatter}

\title{Limit Cycle Bifurcations in~a~Quartic~Ecological~Model\thanksref{label1}}
\thanks[label1]{This work is supported by the Netherlands Organization for Scientific
Research.}

\author{Henk W. Broer}

\address{University of Groningen, Department of Mathematics,
P.O.~Box~407,~9700~AK~Groningen,~The~Netherlands}

\author{Valery A. Gaiko}

\ead{valery.gaiko@yahoo.com}

\address{Belarusian State University of Informatics and Radioelectronics,
Department~of~Mathematics,~L.\,Beda~Str.\,6--4,~Minsk~220040,~Belarus}

\begin{abstract}
In this paper we complete the global qualitative analysis of a quartic ecological model.
In particular, studying global bifurcations of singular points and limit cycles, we prove
that the corresponding dynamical system has at most two limit cycles.
\end{abstract}

\end{frontmatter}

\section{Introduction}

The paper is based on the applications of Bifurcation Theory and can be used for modeling problems, where system parameters play a certain role in various bifurcations. In this paper we consider a particular (quartic) family of planar vector fields, which models the dynamics of the populations of predators and their prey in a given ecological system and which is a variation on the classical Lotka--Volterra system. For the latter system the change of the prey density per unit of time per predator called the response function is proportional to the prey density. This means that there is no saturation of the predator when the amount of available prey is large. However, it is more realistic to consider a nonlinear and bounded response function, and in fact different response functions have been used in the literature to model the predator response, see \cite{Bazykin}--\cite{bnrsw}, \cite{hol}--\cite{lx}, \cite{zcw}.
    \par
For instance, Zhu et al. \cite{zcw} have studied recently the following predator-prey model:
    $$
    \begin{array}{l}
\dot{x}=x(a-\lambda x)-yP(x) \qquad \mbox{(prey)},
    \\[2mm]
\dot{y}=-\delta y + yQ(x) \ \ \ \qquad \mbox{(predator)}.
    \\[2mm]
    \end{array}
    \eqno(1.1)
    $$
The variables $x>0$ and $y>0$ denote the density of the prey and predator populations respectively,
while $P(x)$ is a non-monotonic response function given by
    $$
P(x)=\frac{mx}{\alpha x^{2}+\beta x+1},
    \\[2mm]
    \eqno(1.2)
    $$
where $\alpha,$ $m$ are positive and where $\beta>-2\sqrt{\alpha}.$ Observe that in the absence of predators,
the number of prey increases according to a logistic growth law. The coefficient a represents the
intrinsic growth rate of the prey, while $\lambda>0$ is the rate of competition or resource limitation
of prey. The natural death rate of the predator is given by $\delta>0.$ In Gause's model the function
$Q(x)$ is given by $Q(x)=cP(x),$ where $c>0$ is the rate of conversion between prey and predator.
For further discussion on the biological relevance of the model see \cite{bnrs1}--\cite{bnrsw}.
    \par
In this paper we investigate the following family
    $$
    \begin{array}{l}
\displaystyle\dot{x}=x\left(1-\lambda x-\frac{y}{\alpha x^{2}+\beta x+1}\right) \qquad \qquad \mbox{(prey)},
    \\[4mm]
\displaystyle\dot{y}=y\left(-\delta -\mu y+\frac{x}{\alpha x^{2}+\beta x+1}\right) \qquad \mbox{(predator)},
    \end{array}
    \eqno(1.3)
    $$
where $\alpha\geq0,$ $\delta>0,$ $\lambda>0,$ $\mu\geq0$ and $\beta>-2\sqrt{\alpha}$ are parameters.
We note that (1.3) is obtained from (1.1) by adding the term $-\mu y^{2}$ to the second equation and after scaling
$x$ and $y,$ as well as the parameters and the time $t.$ In this way we take into account competition between predators
for resources other than prey. The non-negative coefficient $\mu$ is the rate of competition amongst predators. For examples
of populations that use the group defense strategy see \cite{Bazykin}--\cite{bnrsw}.
    \par
System (1.3) can be written in the form
    $$
    \begin{array}{l}
\dot{x}=~~x((1-\lambda x)(\alpha x^{2}+\beta x+1)-y)\equiv P,
    \\[2mm]
\dot{y}=-y((\delta +\mu y)(\alpha x^{2}+\beta x+1)-x)\equiv Q.
    \end{array}
    \eqno(1.4)
    $$
This quartic ecological model was studied earlier, for instance, in \cite{bnrs1}--\cite{bnrsw}. However, the qualitative
analysis was incomplete, since the global bifurcations of limit cycles could not be studied properly by means of the
methods and techniques which were used earlier in the qualitative theory of dynamical systems.
    \par
Together with (1.4), we will also consider an auxiliary system (see \cite{BL}, \cite{Perko})
    $$
\dot{x}=P-\gamma Q, \qquad \dot{y}=Q+\gamma P,
    \eqno(1.5)
    $$
applying to these systems new bifurcation methods and geometric approaches developed in \cite{Gaiko}--\cite{gvh}
and completing the qualitative analysis of (1.4).

\section{Preliminaries}

In this paper geometric aspects of Bifurcation Theory are used and de\-ve\-loped. First of all, the two-isocline method which was developed by Erugin is used, see \cite{Gaiko}. An isocline portrait is the most natural construction for a polynomial equation. It is sufficient to have only two isoclines (of zero and infinity) to obtain principal information on the original polynomial system, because these two isoclines are right-hand sides of the system. Geometric properties of isoclines (conics, cubics, quartics, etc.) are well-known, and all isocline portraits can be easily constructed. By means of them, all topologically different qualitative pictures of integral curves to within a number of limit cycles and distinguishing center and focus can be obtained. Thus, it is possible to carry out a rough topological classification of the phase portraits for the polynomial dynamical systems. It is the first application of Erugin's method. After studying contact and rotation properties of the isoclines, the simplest (canonical) systems containing limit cycles can be also constructed. Two groups of parameters can be distinguished in such systems: static and dynamic. Static parameters determine the behavior of phase trajectories in principle, since they control the number, position, and character of singular points in a finite part
of the plane (finite singularities). The parameters from the first group determine also a possible behavior of separatrices and singular points at infinity (infinite singularities) under variation of the parameters from the second group. The dynamic parameters are field rotation parameters, see \cite{BL}, \cite{Gaiko}, \cite{Perko}. They do not change the number, position and index of the finite singularities, but only involve the vector field in a directional rotation. The rotation parameters allow to control the infinite singularities, the behavior of limit cycles and separatrices. The cyclicity of singular points and separatrix cycles, the behavior of semi-stable and other multiple limit cycles are controlled by these parameters as well. Therefore, by means of the rotation parameters,
it is possible to control all limit cycle bifurcations and to solve the most complicated problems of the qualitative theory of dynamical systems.
    \par
In \cite{Gaiko}, \cite{gai1} some complete results on quadratic systems have been presented. In particular, it has been proved that for quadratic systems four is really the maximum number of limit cycles and $(3\!:\!1),$ i.\,e., three limit cycles around one focus and
 the only limit cycle around another focus, is their only possible distribution (this is a solution of Hilbert's Sixteenth Problem in 
 the quadratic case of polynomial dynamical systems). In \cite{gvh} some preliminary results on generalizing new ideas and methods of \cite{Gaiko} to cubic dynamical systems have already been established. In particular, a canonical cubic system of Kukles type has been constructed and the global qualitative analysis of its special case corresponding to a generalized Li\'{e}nard equation has been carried out. It has been proved also that the foci of such a Li\'{e}nard system can be at most of second order and that such system can have at most three limit cycles on the whole phase plane. Moreover, unlike all previous works on the Kukles-type systems, global bifurcations of limit and separatrix cycles using arbitrary (including as large as possible) field rotation parameters of the canonical system have been studied in \cite{gvh}. As a result, the classification of all possible types of separatrix cycles for the generalized Li\'{e}nard system has been obtained and all possible distributions of its limit cycles have been found. In \cite{gai2} a solution of Smale's Thirteenth Problem proving that the Li\'{e}nard system with a polynomial of degree $2k+1$ can have
at most $k$ limit cycles has been presented. All of these methods and results can be applied to quartic dynamical systems as well.
In this paper, using \cite{Gaiko}--\cite{gvh}, we will complete the global qualitative analysis of quartic ecological model (1.4).
In particular, studying global bifurcations of singular points and limit cycles, we will prove that the corresponding dynamical
system has at most two limit cycles.

\section{Singular Points}

The study of singular point of system (1.4) will use two index theorems by H.\,Poincar\'{e}, see \cite{BL}.
But first let us define the Poincar\'{e} index~\cite{BL}.
    \medskip
    \par
    \textbf{Definition 3.1.}
Let $S$ be a simple closed curve in the phase plane not passing
through a singular point of the system
    $$
    \dot{x}=P(x,y), \quad \dot{y}=Q(x,y),
    \eqno(3.1)
    $$
where $P(x,y)$ and $Q(x,y)$ are continuous functions (for example,
polynomials), and $M$ be some point on $S.$ If the point $M$ goes
around the curve $S$ in positive direction (counterclockwise) one
time, then the vector coinciding with the direction of a tangent
to the trajectory passing through the point $M$ is rotated through
the angle $2\pi j$ $(j=0,\pm1,\pm2,\ldots).$ The integer $j$ is
called the \emph{Poincar\'{e} index} of the closed curve $S$
relative to the vector field of system~(3.1) and has the
expression
    $$
    j=\frac{1}{2\pi}\oint_S\frac{P~dQ-Q~dP}{P^2+Q^2}.\\[-2mm]
    $$
    \par
According to this definition, the index of a node or a focus, or a
center is equal to $+1$ and the index of a saddle is $-1.$
    \par
    \medskip
    \textbf{Theorem 3.1 (First Poincar\'{e} Index Theorem).}
    \emph{If $N,$ $N_f,$ $N_c,$ and $C$ are respectively the number
of nodes, foci, centers, and saddles in a finite part of the phase
plane and $N'$ and $C'$ are the number of nodes and saddles at
infinity, then it is valid the formula}
    $$
    N+N_f+N_c+N'=C+C'+1.
    $$
    \par
    \textbf{Theorem 3.2 (Second Poincar\'{e} Index Theorem).}
    \emph{If all singular points are simple, then along an isocline
without multiple points lying in a Poincar\'{e} hemisphere which is 
obtained by a stereographic projection of the phase plane, the
singular points are distributed so that a saddle is followed by a
node or a focus, or a center and vice versa. If two points are
separated by the equator of the Poincar\'{e} sphere, then a saddle
will be followed by a saddle again and a node or a focus, or a
center will be followed by a node or a focus, or a center.}
    \medskip
    \par
We will use also the following theorem by A.\,N.\,Berlinskii, see \cite{Gaiko}.
    \par
    \medskip
    \textbf{Theorem 3.3 (Berlinskii Theorem)}.
\emph{If a quadratic system~$(3.1)$ has four singular points in a finite part
of the phase plane, then only one of the following cases is possible$:$ $
1)\!$~these points are vertices of a convex quadrangular, where
two opposite vertices are saddles $($antisaddles$)$ and two others
are antisaddles $($saddles$);$ $2)\!$~three singular points are
vertices of a triangle containing the fourth point inside, and
if this point is a saddle $($antisaddle$),$ then the others are
antisaddles $($saddles$),$ where antisaddles are singularities 
that are not a saddle.}
    \medskip
    \par
Consider system (1.4) which has two invariant straight lines: $x=0$ and $y=0.$
Its finite singularities are determined by the algebraic system
   $$
    \begin{array}{l}
x((1-\lambda x)(\alpha x^{2}+\beta x+1)-y)=0,
    \\[2mm]
y((\delta+\mu y)(\alpha x^{2}+\beta x+1)-x)=0.\\[2mm]
    \end{array}
    \eqno(3.2)
    $$
From (3.2), we have got: two singular points $(0,0)$ and $(0,-\delta/\mu),$ at most
two points defined by the condition
    $$
\alpha x^{2}+\beta x+1=0, \quad y=0,
    \eqno(3.3)
    $$
and at most four singularities defined by the system
   $$
    \begin{array}{l}
y=(1-\lambda x)(\alpha x^{2}+\beta x+1),
    \\[2mm]
y\,(\delta+\mu y)-x\,(1-\lambda x)=0,
    \end{array}
    \eqno(3.4)
    $$
among which we always have the point $(1/\lambda,0).$
    \par
To investigate the character and distribution of the singular points in the phase plane,
we will use the method developed in \cite{Gaiko}--\cite{gvh}. The sense of this method
is to obtain the simplest (well-known) system by vanishing some parameters (usually field
rotation parameters) of the original system and then to input these parameters successively
one by one studying the dynamics of the singular points (both finite and infinite) in the
phase plane.
    \par
Let the parameters $\alpha,$ $\beta$ vanish and consider first the quadratic system
    $$
    \begin{array}{l}
\dot{x}=~~x(1-\lambda x-y),
    \\[2mm]
\dot{y}=-y(\delta+\mu y-x).
    \end{array}
    \eqno(3.5)
    $$
System (3.5) has four finite singularities, if $\delta\neq1/\lambda.$ Studying isocline portraits
of the equation corresponding to system (3.5) and applying theorems~3.1--3.3, we can see that for
the case, when $\delta>1/\lambda,$ system (3.5) has two saddles: $(0,0)$ and a point of intersection
of two straight line-isoclines:
    $$
1-\lambda x-y=0, \quad \delta+\mu y-x=0,
    \eqno(3.6)
    $$
--- and two nodes: $(0,-\delta/\mu)$ and $(1/\lambda,0).$
For the case, when $\delta<1/\lambda,$ system (3.5) has two saddles: $(0,0)$ and $(1/\lambda,0),$~---
and two nodes: $(0,-\delta/\mu)$ and (3.6). If $\delta=1/\lambda,$ it has three singularities:
a saddle $(0,0),$ a node $(0,-\delta/\mu),$ and a saddle-node $(1/\lambda,0).$ Since we consider
the first coordinate quadrant with respect to the variables $x$ and $y,$ we will be interested
basically in the case of $\delta<1/\lambda,$ when the singular point defined by (3.6) is in the
first quadrant.
    \par
To study singular points at infinity, consider the corresponding differential equation
    $$
\frac{dy}{dx}=-\frac{y(\delta+\mu y-x)}{x(1-\lambda x-y)}.\\[2mm]
    \eqno(3.7)
    $$
Dividing the numerator and denominator of the right-hand side of (3.7) by $x^{2}$ $(x\neq0)$ and
denoting $y/x$ by $u,$ we will get the algebraic equation
    $$
(1-\mu)u^{2}+(1+\lambda)u=0, \quad \mbox{where} \quad u=y/x,
    \eqno(3.8)
    $$
for all infinite singularities of (3.7) except when $x=0$ (the ``ends'' of the $y$-axis), see \cite{BL},
\cite{Gaiko}. For this special case we can divide the numerator and denominator of the right-hand side of (3.7)
by $y^{2}$ $(y\neq0)$ denoting $x/y$ by $v$ and consider the algebraic equation
    $$
(1+\lambda)v^{2}+(1-\mu)v=0, \quad \mbox{where} \quad v=x/y.
    \eqno(3.9)
    $$
The equations (3.8) and (3.9) give three singular points at infinity for (3.7): two nodes on
the ``ends'' of the $x$ and $y$ axes and a saddle in the direction of $u=(\lambda+1)/(\mu-1).$
    \par
Fix the parameters $\delta,$ $\lambda,$ $\mu$ and take $\beta<0$ (this case will be also more
inte\-resting to us). After inputting the parameter $\beta,$ we will have a cubic system:
    $$
    \begin{array}{l}
\dot{x}=~~x((1-\lambda x)(\beta x+1)-y),
    \\[2mm]
\dot{y}=-y((\delta+\mu y)(\beta x+1)-x).
    \end{array}
    \eqno(3.10)
    $$
For $\delta<1/\lambda$ and $\beta<0,$ system (3.10) has five finite singularities: two saddles~---
$(0,0)$ and $(1/\lambda,0),$ two nodes~--- $(0,-\delta/\mu)$ and $(-1/\beta,0),$ and an antisaddle
(a node, a focus, or a center) defined as a point of intersection of two isoclines:
    $$
    \begin{array}{l}
(1-\lambda x)(\beta x+1)-y=0,\\[1mm]
(\delta+\mu y)(\beta x+1)-x=0.\\[2mm]
    \end{array}
    \eqno(3.11)
    $$
For singular points at infinity, consider the corresponding differential equation
    $$
\frac{dy}{dx}=-\frac{y((\delta+\mu y)(\beta x+1)-x)}{x((1-\lambda x)(\beta x+1)-y)}
    \eqno(3.12)
    $$
and the algebraic equations
    $$
\mu u^{2}-\lambda u=0, \quad \mbox{where} \quad u=y/x,\\[-2mm]
    \eqno(3.13)
    $$
and
    $$
\lambda v^{3}-\mu v^{2}=0, \quad \mbox{where} \quad v=x/y,\\[2mm]
    \eqno(3.14)
    $$
which give three infinite singularities: a node on the ``ends'' of the $x$-axis, a saddle-node
on the ``ends'' of the $y$-axis, and a saddle in the direction of $u=\lambda/\mu.$
    \par
Fix the parameters $\beta,$ $\delta,$ $\lambda,$ $\mu$ and take $\alpha>0.$ Studying the bundle of
cubic curves
    $$
y=(1-\lambda x)(\alpha x^{2}+\beta x+1)\\[2mm]
    \eqno(3.15)
    $$
which intersect at the point $(1/\lambda,0)$ and contact at the point $(0,1),$ we can see that
system (1.4) obtained after inputting $\alpha$ has first six finite singular points: three saddles~---
$(0,0),$ $(1/\lambda,0),$ and $\displaystyle((-\beta+\sqrt{\beta^{2}-4\alpha}\,)/(2\alpha),0),$
two nodes~--- $(0,-\delta/\mu)$ and $\displaystyle((-\beta-\sqrt{\beta^{2}-4\alpha}\,)/(2\alpha),0),$
and an antisaddle defined as a point of intersection of isoclines (3.4).
    \par
On increasing the parameter $\alpha,$ the points $\displaystyle((-\beta-\sqrt{\beta^{2}-4\alpha}\,)/(2\alpha),0)$
and $\displaystyle((-\beta+\sqrt{\beta^{2}-4\alpha}\,)/(2\alpha),0)$ combine a saddle-node which then disappears.
On further increasing $\alpha,$ the point $(1/\lambda,0)$ becomes a triple saddle from which a saddle
and a node (or a saddle-node) will appear. Thus, we will have three singular points in the first quadrant:
a saddle $S$ and two antisaddles~--- $A_{1}$ and $A_{2}$ which are defined as points of intersection of isoclines
(3.4). Suppose that with respect to the $x$-axis they have the following sequence: $A_{1},$ $S,$ $A_{2}.$
    \par
To study singular points of (1.4) at infinity, consider the corresponding diffe\-rential equation
    $$
\frac{dy}{dx}=-\frac{y((\delta+\mu y)(\alpha x^{2}+\beta x+1)-x)}{x((1-\lambda x)(\alpha x^{2}+\beta x+1)-y)}\\[2mm]
    \eqno(3.16)
    $$
and the algebraic equations
    $$
\mu u^{2}-\lambda u=0, \quad \mbox{where} \quad u=y/x,\\[-2mm]
    \eqno(3.17)
    $$
and
    $$
\lambda v^{4}-\mu v^{3}=0, \quad \mbox{where} \quad v=x/y,\\[2mm]
    \eqno(3.18)
    $$
which give three infinite singularities: a simple node on the ``ends'' of the\linebreak $x$-axis, a triple node
on the ``ends'' of the $y$-axis, and a simple saddle in the direction of $u=\lambda/\mu.$
    \par
Note that all results on finite singularities of system (1.4) agree with the results of \cite{bnrs1}--\cite{bnrsw},
\cite{lx}, \cite{zcw}, but where infinite singularities have not been investigated at all. Using the obtained information
and applying the approach developed in \cite{Gaiko}--\cite{gvh}, we can study limit cycle bifurcations of system (1.4) now.
This study will use also some results obtained in \cite{bnrs1}--\cite{bnrsw}, \cite{lx}, \cite{zcw}. In particular, the results
on the cyclicity of singular points of (1.4) will be used. However, it is surely not enough to have only these results to prove
the main theorem of this paper: on the maximum number of limit cycles of (1.4).

\section{Bifurcations of Limit Cycles}

\noindent Let us first formulate the Wintner--Perko termination principle~\cite{Perko}
for the polynomial system
    $$
    \mbox{\boldmath$\dot{x}$}=\mbox{\boldmath$f$}
    (\mbox{\boldmath$x$},\mbox{\boldmath$\mu$)},\\[2mm]
    \eqno(4.1_{\mbox{\boldmath$\mu$}})
    $$
where $\mbox{\boldmath$x$}\in\textbf{R}^2;$ \
$\mbox{\boldmath$\mu$}\in\textbf{R}^n;$ \
$\mbox{\boldmath$f$}\in\textbf{R}^2$ \ $(\mbox{\boldmath$f$}$ is a
polynomial vector function).
    \par
    \medskip
\noindent\textbf{Theorem 4.1 (Wintner--Perko termination principle).}
    \emph{Any one-para\-me\-ter fa\-mi\-ly of multip\-li\-city-$m$
limit cycles of relatively prime polynomial system
$(4.1_{\mbox{\boldmath$\mu$}})$ can be extended in a unique way to
a maximal one-parameter family of multiplicity-$m$ limit cycles of
$(4.1_{\mbox{\boldmath$\mu$}})$ which is either open or cyclic.}
    \par
\emph{If it is open, then it terminates either as the parameter or
the limit cycles become unbounded; or, the family terminates
either at a singular point of $(4.1_{\mbox{\boldmath$\mu$}}),$
which is typically a fine focus of multiplicity~$m,$ or on a
$($compound$\,)$ separatrix cycle of
$(4.1_{\mbox{\boldmath$\mu$}}),$ which is also typically of
multiplicity~$m.$}
    \medskip
    \par
The proof of this principle for general polynomial system
$(4.1_{\mbox{\boldmath$\mu$}})$ with a vector parameter
$\mbox{\boldmath$\mu$}\in\textbf{R}^n$ parallels the proof of the
pla\-nar termination principle for the system
    $$
    \vspace{1mm}
    \dot{x}=P(x,y,\lambda),
        \quad
    \dot{y}=Q(x,y,\lambda)\\[2mm]
    \eqno(4.1_{\lambda})
    $$
with a single parameter $\lambda\in\textbf{R}$ (see \cite{Gaiko},
\cite{Perko}), since there is no loss of generality in assuming
that system $(4.1_{\mbox{\boldmath$\mu$}})$ is parameterized by a
single parameter $\lambda;$ i.\,e., we can assume that there
exists an analytic mapping $\mbox{\boldmath$\mu$}(\lambda)$ of
$\textbf{R}$ into $\textbf{R}^n$ such that
$(4.1_{\mbox{\boldmath$\mu$}})$ can be written as
$(4.1\,_{\mbox{\boldmath$\mu$}(\lambda)})$ or even
$(4.1_{\lambda})$ and then we can repeat everything, what had been
done for system $(4.1_{\lambda})$ in~\cite{Perko}. In particular,
if $\lambda$ is a field rotation parameter of $(4.1_{\lambda}),$
the following Perko's theorem on monotonic families of limit cycles
is valid (see \cite{Perko}).
    \par
    \medskip
\noindent\textbf{Theorem 4.2.}
    \emph{If $L_{0}$ is a nonsingular multiple limit cycle of
$(4.1_{0}),$ then  $L_{0}$ belongs to a one-parameter family of
limit cycles of $(4.1_{\lambda});$ furthermore\/$:$}
    \par
1)~\emph{if the multiplicity of $L_{0}$ is odd, then the family
either expands or contracts mo\-no\-to\-ni\-cal\-ly as $\lambda$
increases through $\lambda_{0};$}
    \par
2)~\emph{if the multiplicity of $L_{0}$ is even, then $L_{0}$
bi\-fur\-cates into a stable and an unstable limit cycle as
$\lambda$ varies from $\lambda_{0}$ in one sense and $L_{0}$
dis\-ap\-pears as $\lambda$ varies from $\lambda_{0}$ in the
opposite sense; i.\,e., there is a fold bifurcation at
$\lambda_{0}.$}
    \medskip
    \par
Applying the definition of a field rotation parameter  \cite{BL}, \cite{Gaiko}, \cite{Perko},
i.\,e., a parameter which rotates the field in one direction, to system (1.4), let us calculate
the corresponding determinants for the parameters $\alpha$ and $\beta,$ respectively:
    $$
\Delta_{\alpha}=PQ'_{\alpha}-QP'_{\alpha}=x^{3}y(y(\delta+\mu y)-x(1-\lambda x)),
    \eqno(4.2)
    $$
    $$
\Delta_{\beta}=PQ'_{\beta}-QP'_{\beta}=x^{2}y(y(\delta+\mu y)-x(1-\lambda x)).\\[4mm]
    \eqno(4.3)
    $$
It follows from (4.2) and (4.3) that on increasing $\alpha$ or $\beta$ the vector field of (1.4) in the
first quadrant is rotated in positive direction (counterclockwise) only on the outside of the ellipse
    $$
y(\delta+\mu y)-x(1-\lambda x)=0.
    \eqno(4.4)
    $$
Therefore, to study limit cycle bifurcations of system (1.4), it makes sense together with (1.4) to consider
also an auxiliary system (1.5) with a field rotation parameter $\gamma\!:$
    $$
\Delta_{\gamma}=P^{2}+Q^{2}\geq0.
    \eqno(4.5)
    $$
    \par
Using system (1.5) and applying Perko's results, we will prove the following theorem.
	\par
    \medskip
\noindent \textbf{Theorem 4.3.}
\emph{System~$(1.4)$ has at most two limit cycles.}
	\par
    \medskip
\noindent\textbf{Proof.} First let us prove that system (1.4) can have at least two limit cycles.
Begin with quadratic system (3.5). It is clear that such a system, with two invariant straight lines,
cannot have limit cycles at all \cite{Gaiko}. Inputting a negative parameter $\beta$ into this system,
the vector field of cubic system (3.10) will be rotated in negative direction (clockwise) at infinity, the structure
and the character of stability of infinite singularities will be changed, and an unstable limit, $\Gamma_{1},$ will
appear immediately from infinity in this case. This cycle will surround a stable antisaddle (a node or a focus) $A_{1}$
which is in the first quadrant of system (3.10). Inputting a positive parameter $\alpha$ into system (3.10),
the vector field of quartic system (1.4) will be rotated in positive direction (counterclockwise) at infinity,
the structure and the character of stability of infinite singularities will be changed again, and a stable limit,
$\Gamma_{2},$ surrounding $\Gamma_{1}$ will appear immediately from infinity in this case. On further increasing
the parameter $\alpha,$ the limit cycles $\Gamma_{1}$ and $\Gamma_{2}$ combine a semi-stable limit, $\Gamma_{12},$
which then disappears in a ``trajectory concentration'' \cite{BL}, \cite{Gaiko}.
    \par
As we saw above, on further increasing $\alpha,$ two other singular points, a saddle $S$ and an antisaddle $A_{2},$
will appear in the first quadrant in system (1.4). We can fix the parameter $\alpha,$ fixing simultaneously the positions
of the finite singularities $A_{1},$ $S,$ $A_{2},$ and consider system (1.5) with a positive parameter~$\gamma$ which
acts like a positive parameter $\alpha$ of system (1.4), but on the whole phase plane.
    \par
So, consider system (1.5) with a positive parameter~$\gamma.$ On increasing this parameter, the stable nodes $A_{1}$
and $A_{2}$ becomes first stable foci, then they change the character of their stability, becoming unstable foci.
At these Andronov--Hopf bifurcations \cite{BL}, \cite{Gaiko}, stable limit cycles will appear from the foci $A_{1}$ and
$A_{2}.$ On further increasing~$\gamma,$ the limit cycles will expand and will disappear in small separatrix loops of
the saddle $S.$ If these loops are formed simultaneously, we will have a so-called eight-loop separatrix cycle.
In this case, a big stable limit surrounding three singular points, $A_{1},$ $S,$ and $A_{2},$ will appear
from the eight-loop separatrix cycle after its destruction, expanding to infinity on increasing~$\gamma.$
If a small loop is formed earlier, for example, around the point $A_{1}$ $(A_{2}),$ then, on increasing~$\gamma,$
a big loop formed by two lower (upper) adjoining separatrices of the saddle~$S$ and surrounding the points $A_{1}$
and $A_{2}$ will appear. After its destruction, we will have simultaneously a big limit cycle surrounding
three singular points, $A_{1},$ $S,$ $A_{2},$ and a small limit cycle surrounding the point $A_{2}$ $(A_{1}).$
Thus, we have proved that system (1.4) can have at least two limit cycles, see also \cite{bnrs1}--\cite{bnrsw},
\cite{lx}, \cite{zcw}.
    \par
Let us prove now that this system has at most two limit cycles. The proof is carried out by contradiction applying
Catastrophe Theory, see \cite{Gaiko}, \cite{Perko}. Consider system (1.5) with three parameters: $\alpha,$ $\beta,$
and $\gamma$ (the parameters $\delta,$ $\lambda,$ and $\mu$ can be fixed, since they do not generate limit cycles).
Suppose that (1.5) has three limit cycles surrounding the only point, $A_{1},$ in the first quadrant.
Then we get into some domain of the parameters $\alpha,$ $\beta,$ and $\gamma$ being restricted
by definite con\-di\-tions on three other parameters, $\delta,$ $\lambda,$ and $\mu$. This domain
is bounded by two fold bifurcation surfaces forming a cusp bifurcation surface of multiplicity-three
limit cycles in the space of the pa\-ra\-me\-ters $\alpha,$ $\beta,$ and $\gamma$ \cite{Gaiko}, \cite{Perko}.
    \par
The cor\-res\-pon\-ding maximal one-parameter family of multiplicity-three limit cycles cannot be cyclic, otherwise there
will be at least one point cor\-res\-pon\-ding to the limit cycle of multi\-pli\-ci\-ty four (or even higher) in the parameter
space. Extending the bifurcation curve of multi\-pli\-ci\-ty-four limit cycles through this point and parameterizing the
corresponding maximal one-parameter family of multi\-pli\-ci\-ty-four limit cycles by the field rotation para\-me\-ter,
$\gamma,$ according to Theorem~4.2, we will obtain two monotonic curves of multi\-pli\-ci\-ty-three and one, respectively,
which, by the Wintner--Perko termination principle (Theorem~4.1), terminate either at the point $A_{1}$ or on a separatrix
cycle surrounding this point. Since we know at least the cyclicity of the singular point which is equal to two (see \cite{bnrs1}--\cite{bnrsw}, \cite{lx}, \cite{zcw}), we have got a contradiction with the termination principle stating
that the multiplicity of limit cycles cannot be higher than the multi\-pli\-ci\-ty (cyclicity) of the singular point
in which they terminate.
    \par
If the maximal one-parameter family of multiplicity-three limit cycles is not cyclic, using the same principle (Theorem~4.1),
this again contradicts the cyclicity of $A_{1}$ (see \cite{bnrs1}--\cite{bnrsw}, \cite{lx}, \cite{zcw}) not admitting the
multiplicity of limit cycles to be higher than two. This contradiction completes the proof in the case of one singular
point in the first quadrant.
    \par
Suppose that system (1.5) with three finite singularities, $A_{1},$ $S,$ and $A_{2},$ has two small limit cycles around,
for example, the point $A_{1}$ (the case when limit cycles surround the point $A_{2}$ is considered in a similar way).
Then we get into some domain in the space of the parameters $\alpha,$ $\beta,$ and $\gamma$ which is bounded by a fold
bifurcation surface of multiplicity-two limit cycles \cite{Gaiko}, \cite{Perko}.
    \par
The cor\-res\-pon\-ding maximal one-parameter family of multiplicity-two limit cycles cannot be cyclic, otherwise there
will be at least one point cor\-res\-pon\-ding to the limit cycle of multi\-pli\-ci\-ty three (or even higher) in the
parameter space. Extending the bifurcation curve of multi\-pli\-ci\-ty-three limit cycles through this point and
parameterizing the corresponding maximal one-parameter family of multi\-pli\-ci\-ty-three limit cycles by the field
rotation para\-me\-ter, $\gamma,$ according to Theorem~4.2, we will obtain a monotonic curve which, by the Wintner--Perko
termination principle (Theorem~4.1), terminates either at the point $A_{1}$ or on some separatrix cycle surrounding this
point. Since we know at least the cyclicity of the singular point which is equal to one in this case \cite{bnrs1}--\cite{bnrsw},
\cite{lx}, \cite{zcw}, we have got a contradiction with the termination principle (Theorem~4.1).
    \par
If the maximal one-parameter family of multiplicity-two limit cycles is not cyclic, using the same principle (Theorem~4.1),
this again contradicts the cyclicity of $A_{1}$ (see \cite{bnrs1}--\cite{bnrsw}, \cite{lx}, \cite{zcw}) not admitting the
multiplicity of limit cycles higher than one. Moreover, it also follows from the termination principle that either an ordinary
(small) separatrix loop or a big loop, or an eight-loop cannot have the multiplicity (cyclicity) higher than one in this case.
Therefore, according to the same principle, there are no more than one limit cycle in the exterior domain surrounding all
three finite singularities, $A_{1},$ $S,$ and $A_{2}.$
    \par
Thus, taking into account all other possibilities for limit cycle bifurcations (see \cite{bnrs1}--\cite{bnrsw},
\cite{lx}, \cite{zcw}), we conclude that system~(1.4) cannot have either a multiplicity-three limit cycle or
more than two limit cycles in any configuration. The theorem is proved.
\qquad $\Box$

\end{document}